\documentclass{amsart}
\usepackage{graphicx}
\usepackage{amscd}             
\usepackage[curve]{xypic}

\newbox\mybox
\def\overtag#1#2#3{\setbox\mybox\hbox{$#1$}\hbox to
  0pt{\vbox to 0pt{\vglue-#3\vglue-\ht\mybox\hbox to \wd\mybox
      {\hss$\ss#2$\hss}\vss}\hss}\box\mybox}
\def\undertag#1#2#3{\setbox\mybox\hbox{$#1$}\hbox to 0pt{\vbox to
    0pt{\vglue#3\vglue\ht\mybox\hbox to \wd\mybox
      {\hss$\ss#2$\hss}\vss}\hss}\box\mybox}
\def\lefttag#1#2#3{\hbox to 0pt{\vbox to 0pt{\vss\hbox to
      0pt{\hss$\ss#2$\hskip#3}\vss}}#1}
\def\righttag#1#2#3{\hbox to 0pt{\vbox to 0pt{\vss\hbox to
      0pt{\hskip#3$\ss#2$\hss}\vss}}#1}
\let\ss\scriptstyle

\def\Dot{\lower.2pc\hbox to 2.5pt{\hss$\bullet$\hss}}
\def\Circ{\lower.2pc\hbox to 2.5pt{\hss$\circ$\hss}}
\def\Vdots{\raise5pt\hbox{$\vdots$}}
\def\splicediag#1#2{\xymatrix@R=#1pt@C=#2pt@M=0pt@W=0pt@H=0pt}
\newcommand\lineto{\ar@{-}}
\newcommand\dashto{\ar@{--}}
\newcommand\dotto{\ar@{.}}

\usepackage{pinlabel}

\newcommand\bicolor{two-color}
\newcommand\Bicolor{Two-color}

\renewcommand{\H}{{\mathbb H}}
\newcommand{\Z}{{\mathbb Z}}
\newcommand{\R}{{\mathbb R}}

\newcommand{\E}{{\mathbb E}}
\renewcommand{\S}{{\mathbb S}}

\newcommand{\GG}{\mathcal {G}}

\newcommand{\co}{\colon}

\usepackage{color}

\newtheorem{theorem}{Theorem}[section]
\newtheorem*{theorem*}{Theorem}

\newtheorem{lemma}[theorem]{Lemma}
\newtheorem{proposition}[theorem]{Proposition}

\theoremstyle{definition}
\newtheorem{definition}[theorem]{Definition}

\newtheorem{remark}[theorem]{Remark}

\newtheorem{question}[theorem]{Question}

\long\def\Restate#1#2#3#4{
\medskip\par\noindent
{\bf #1 \ref{#2}{\rm #3}.}
{\it #4}\par\medskip
}

\begin{document}
\title
[Quasi-isometric classification of graph manifolds] {Quasi-isometric
  classification of graph manifold groups} 
\author%[Behrstock]
{Jason A. Behrstock}
\address{Department of Mathematics\\The University of Utah\\Salt Lake
  City, UT 84112, USA} 
\email{jason@math.utah.edu}
\author%[Neumann]
{Walter D.
  Neumann} \thanks{Research supported under NSF grants no.\
  DMS-0083097 and DMS-0604524} 
  \address{Department of Mathematics\\Barnard College,
  Columbia University\\New York, NY 10027, USA}
\email{neumann@math.columbia.edu} 
\begin{abstract}
  We show that the fundamental groups of any two closed irreducible
  non-geometric graph-manifolds are quasi-isometric. This answers a
  question of Kapovich and Leeb. We also classify the quasi-isometry
  types of fundamental groups of graph-manifolds with boundary in terms
  of certain finite two-colored graphs.  A corollary is the 
  quasi-isometric classification of Artin groups whose presentation
  graphs are trees.  In particular any two right-angled Artin groups
  whose presentation graphs are trees of diameter greater than $2$ are
  quasi-isometric, answering a question of Bestvina;  further, this
  quasi-isometry class does not include any other right-angled Artin
  groups.
\end{abstract}
\maketitle

A finitely generated group can be considered geometrically when 
endowed with a word metric---up to quasi-isometric equivalence, 
such metrics are unique. (Henceforth only finitely generated 
groups will be considered.) Given a collection of groups,~$\GG$, 
Gromov proposed the fundamental questions of identifying which  
groups are quasi-isometric to those in $\GG$ 
(\emph{rigidity}) and which groups in $\GG$ are 
quasi-isometric to each other (\emph{classification}) 
\cite{Gromov:Asymptotic}.

In this paper, we focus on the classification question for graph
manifold groups and right-angled Artin groups.

A compact $3$-manifold $M$ is called \emph{geometric} if $M\setminus 
\partial M$ admits a geometric structure in the sense of Thurston
(i.e., a complete locally homogeneous Riemannian metric of finite
volume).  The Geometrization Conjecture \cite{Perelman:Geom1,
  Perelman:Geom2, Perelman:Geom3} provides that every irreducible
$3$-manifold of zero Euler characteristic (i.e., with boundary
consisting only of tori and Klein bottles) admits a decomposition
along tori and Klein bottles into geometric pieces, the minimal such
decomposition is called the \emph{geometric decomposition}.

There is a considerable literature on quasi-isometric rigidity and
classification of 3-manifold groups.  The rigidity results can be
briefly summarized:

\begin{theorem*}
  If a group  $G$ is quasi-isometric to the fundamental group
  of a 3-manifold $M$ with zero Euler characteristic, then $G$ is
  weakly commensurable\footnote{Two groups are said to be \emph{weakly
      commensurable} if they have quotients by finite normal subgroups
    which have isomorphic finite index subgroups.}  with $\pi_1(M')$
  for some such 3-manifold $M'$. Moreover, $M'$ is closed resp.\
  irreducible resp.\ geometric if and only if $M$ is.
\end{theorem*}

This quasi-isometric rigidity for $3$-manifold groups is the
culmination of the work of many authors, key steps being provided by
Gromov-Sullivan, Cannon-Cooper, Eskin-Fisher-Whyte, Kapovich-Leeb,
Rieffel, Schwartz \cite{CannonCooper, EskinFisherWhyte,
Gromov:PolynomialGrowth, KapovichLeeb:haken,
Rieffel:H2crossR, Schwartz:RankOne}. The reducible case
reduces to the irreducible case using Papasoglu and Whyte
\cite{PapasogluWhyte:ends} and the irreducible non-geometric case is
considered by Kapovich and Leeb~\cite{KapovichLeeb:haken}.

The classification results in the geometric case 
can be summarized by the following; the first half is an easy
application of the Milnor-\v Svarc Lemma \cite{milnor}, \cite{svarc}:

\begin{theorem*}
  There are exactly seven quasi-isometry classes of fundamental groups
  of closed geometric 3-manifolds, namely any such group is
  quasi-isometric to one of the eight Thurston geometries ($\S^3$,
  $\S^2\times \E^1$, $\E^3$, $\mathrm{Nil}$, $\H^2\times \E^1$,
  $\widetilde{\mathrm{PSL}}$, $\mathrm{Sol}$, $\H^3$) but the two
  geometries $\H^2\times \E^1$ and  $\widetilde{\mathrm{PSL}}$ are
  quasi-isometric.

  If a geometric manifold $M$ has boundary, then it is either Seifert
  fibered and its fundamental group is quasi-isometric (indeed
  commensurable) with $F_2\times \Z$ \cite{KapovichLeeb:haken}, or it
  is hyperbolic, in which case quasi-isometry also implies
  commensurability~\cite{Schwartz:RankOne}.
\end{theorem*}

A \emph{graph manifold} is a $3$-manifold that can be decomposed along
embedded tori and Klein bottles into finitely many Seifert manifolds;
equivalently, these are exactly the class of manifolds with no
hyperbolic pieces in their geometric recomposition.  Since the presence
of a hyperbolic piece can be quasi-isometrically detected
\cite{Gersten:divergence3mflds} \cite{KapovichLeeb:3manifolds}
\cite{BehrstockDrutuMosher:thick}, this implies that the class of
fundamental groups of graph manifolds is rigid. We answer the
classification question for graph manifold groups. 
Before discussing the general case we note the answer for 
closed non-geometric graph manifolds,
resolving a question of Kapovich and Leeb
\cite{KapovichLeeb:3manifolds}.  
\Restate{Theorem} {graphsqi}{}{Any
  two closed non-geometric graph manifolds have bilipschitz
  homeomorphic universal covers. In particular, their fundamental
  groups are quasi-isometric.}

This contrasts with commensurability of closed graph manifolds:
already in the case that the graph manifold is composed of just two
Seifert pieces there are infinitely many commensurability classes
(they are classified in that case but not in general, see Neumann
\cite{Neumann:commens}).

We also classify compact graph manifolds with boundary.  To describe
this we need some terminology. We associate to the geometric
decomposition of a non-geometric graph manifold $M$ its
\emph{decomposition graph} $\Gamma(M)$ which has a vertex for each
Seifert piece and an edge for each decomposing torus or Klein bottle.
We color the vertices of $\Gamma(M)$ \underline{\bf b}lack or
\underline{\bf w}hite according to whether the Seifert piece includes
a boundary component of $M$ or not (\underline{\bf b}ounded or
\underline{\bf w}ithout boundary). We call this the \emph{\bicolor{}ed
  decomposition graph}. We can similarly associate a \bicolor{}ed tree to
the decomposition of the universal cover $\tilde M$ into its fibered
pieces. We call this infinite valence \bicolor{}ed tree $BS(M)$, since
it is the Bass-Serre tree corresponding to the graph of groups
decomposition of $\pi_1(M)$.

The Bass-Serre tree $BS(M)$ can be constructed directly from the
decomposition graph $\Gamma=\Gamma(M)$ by first replacing each edge of
$\Gamma$ by a countable infinity of edges with the same endpoints and
then taking the universal cover of the result. If two \bicolor{}ed
graphs $\Gamma_1$ and $\Gamma_2$ lead to isomorphic \bicolor{}ed trees
by this procedure we say $\Gamma_1$ and $\Gamma_2$ are
\emph{bisimilar}. In Section \ref{sec:bicolor} we give a simpler,
algorithmically checkable, 
criterion for bisimilarity\footnote{We thank Ken Shan for pointing out 
that our equivalence relation is a special case of the computer 
science concept bisimilarity, related to bisimulation.} and show that
each bisimilarity class contains a unique minimal element.

Our classification theorem, which includes the closed case (Theorem
\ref{graphsqi}), is:

\Restate{Theorem}{th:main}{}
  {If $M$ and $M'$ are non-geometric graph manifolds then the
  following are equivalent:
  \begin{enumerate}
  \item $\tilde M$ and $\tilde M'$ are bilipschitz homeomorphic.
  \item $\pi_1(M)$ and $\pi_1(M')$ are quasi-isometric.
\item $BS(M)$ and $BS(M')$ are isomorphic as \bicolor{}ed trees.
\item The minimal \bicolor{}ed graphs in the bisimilarity classes of 
    the decomposition graphs $\Gamma(M)$ and $\Gamma(M')$ are
    isomorphic.
  \end{enumerate}}

One can list minimal \bicolor{}ed graphs of small size,
yielding, for instance, that there are exactly
$2,6,26,199,2811,69711,2921251,204535126,\dots$ quasi-isometry classes of
fundamental groups of non-geometric graph manifolds that are composed
of at most $1,2,3,4,5,6,7,8,\dots$ Seifert pieces.

For closed non-geometric graph manifolds we recover that there is just
one quasi-isometry class (Theorem \ref{graphsqi}): the minimal
\bicolor{}ed graph is a single white vertex with a loop.  Similarly,
for non-geometric graph manifolds that have boundary components in
every Seifert component there is just one quasi-isometry class (the
minimal \bicolor{}ed graph is a single black vertex with a loop).

For graph manifolds with boundary the commensurability classification
is also rich, but not yet well understood. If $M$ consists of two
Seifert components glued to each other such that $M$ has boundary
components in both Seifert components one can show that $M$ is
commensurable with any other such $M$, but this appears to be already
no longer true in the case of three Seifert components.

We end by giving an application to the quasi-isometric classification
of Artin groups. The point is that if the presentation graph is a tree
then the group is a graph-manifold group, so our results apply. In
particular, we obtain the classification of right-angled Artin groups
whose presentation graph is a tree, answering a question of Bestvina.
We also show rigidity of such groups amongst right angled Artin groups.

We call a right-angled Artin group whose presentation graph is a tree a
\emph{right-angled tree group}. If the tree has diameter $\le 2$ then
the group is $\Z,\Z^2$ or (free)$\times\Z$.  We answer Bestvina's
question by showing that right-angled tree groups with presentation
graph of diameter $>2$ are all quasi-isometric to each other. In fact:
\Restate{Theorem}{treegroupsclassification}{} {Let $G'$ be any Artin
  group and let $G$ be a right-angled tree group whose tree has
  diameter $>2$.  Then $G'$ is quasi-isometric to $G$ if and only if
  $G'$ has presentation graph an even-labeled tree of diameter $\ge2$
  satisfying: (i) all interior edges have label 2; and (ii) if the
  diameter is $2$ then at least one edge has label $>2$. (An ``interior
  edge'' is an edge that does not end in a leaf of the tree.)}

The commensurability classification of right-angled tree groups is
richer: Any two whose presentation graphs have diameter~$3$ are
commensurable, but it appears that there are already infinitely many
commensurability classes for diameter~$4$.
  
Theorem~\ref{treegroupsclassification} also has implications for
quasi-isometric rigidity phenomena in relatively hyperbolic groups.
For such applications see Behrstock--Dru\c{t}u--Mosher
\cite{BehrstockDrutuMosher:thick}, where it is shown that graph
manifolds, and thus tree groups, can only quasi-isometrically embed in
relatively hyperbolic groups in very constrained ways.

In the course of proving Theorem~\ref{treegroupsclassification} 
we classify which Artin groups are quasi-isometric 
to $3$-manifold groups. This family of groups coincides with those 
proven by Gordon to be isomorphic to $3$-manifold groups 
\cite{Gordon:Coherence}. 
\subsection*{Acknowledgments} We express our appreciation
to Mladen Bestvina, Mohamad Hindawi, Tadeusz Januszkiewicz and Misha
Kapovich for useful conversations and the referee for pertinent comments.
\section{Quasi-isometry of fattened trees}\label{sec:fattrees}

Let $T$ be a tree all of whose vertices have valence in the interval
$[3, K]$ for some $K$. We fix a positive constant $L$, and 
assume $T$ has been given a simplicial
metric in which each edge has length between $1$ and $L$.  Now
consider a ``fattening'' of $T$, where we replace each edge $E$ by a
strip isometric to $E\times [-\epsilon,\epsilon]$ for some
$\epsilon>0$ and replace each vertex by a regular polygon around the
boundary of which the strips of incoming edges are attached in some
order.  Call this object $X$. Let $X_0$ be similarly constructed, but
starting from the regular 3-valence tree with all edges having length
1, and with $\epsilon=1/2$.

We first note the easy lemma:
\begin{lemma}
There exists $C$, depending only on $K, L,\epsilon$, such that $X$
is $C$-bilipschitz homeomorphic to $X_0$.
  \qed
\end{lemma}

Note that if $S$ is a compact riemannian surface with boundary having
Euler characteristic $<0$ then its universal cover $\tilde S$ is
bilipschitz homeomorphic to a fattened tree as above, and hence to
$X_0$. We can thus use $X_0$ as a convenient bilipschitz model for any
such $\tilde S$.

Let $X$ be a manifold as above, bilipschitz equivalent to $X_0$ (so
$X$ may be a fattened tree or an $\tilde S$). $X$ is a $2$-manifold
with boundary, and its boundary consists of infinitely many copies of
$\R$.
\begin{theorem}\label{fattrees}
  Let $X$ be as above with a chosen boundary component $\partial_0X$.
  Then there exists $K$ and a function $\phi\colon \R\to \R$ such that
  for any $K'$ and any $K'$-bilipschitz homeomorphism $\Phi_0$ from
  $\partial_0X$ to a boundary component $\partial_0X_0$ of the
  ``standard model'' $X_0$, $\Phi_0$ extends to a
  $\phi(K')$-bilipschitz homeomorphism $\Phi\colon X\to X_0$ which is
  $K$-bilipschitz on every other boundary component.
\end{theorem}
\begin{proof}
  If true for some $X$, then the theorem will be true (with $K$
  replaced by $KL$) for any $X'$ $L$-bilipschitz homeomorphic to $X$,
  so we may assume $X$ is isometric to our standard model $X_0$.  In
  this case we will see that $K$ can be arbitrarily close to~$1$ (with
  very slightly more effort one can make $K=1$).

  We will construct the homeomorphism in two steps. The first step
  will be to extend near $\partial_0X$ and the second to extend over
  the rest of $X$. 

  We consider vertices of the underlying tree adjacent to the boundary
  component $\partial_0X$. These will have a certain ``local density''
  along $\partial_0X$ given by the number of them in an interval of a
  given length, measured with respect to the metric on $\partial_0X$
  that pulls back from $\partial_0X_0$ by $\Phi_0$. We first describe
  how to modify this local density using a $(1,L)$-quasi-isometric
  bilipschitz homeomorphism with $L=O(|\log(D)|)$, where $D$ is the
  factor by which we want to modify density.  We increase density
  locally by moves on the underlying tree in which we take a vertex
  along $\partial_0X$ and a vertex adjacent to it not along
  $\partial_0X$ and collapse the edge between them to give a vertex of
  valence 4, which we then expand again to two vertices of valence 3,
  now both along $\partial_0X$ (see Fig.~\ref{Figure 1}).  This can be
  realized by a piecewise-linear homeomorphism. Since it is an
  isometry outside a bounded set, it has a finite bilipshitz bound $k$
  say.  
  To increase density along an interval by at
  most a factor of $D$ we need to repeat this process at most
  $\log_2(D)$ times, so we get a bilipshitz homeomorphism whose
  bilipschitz bound is bounded in terms of $D$. Similarly we can
  decrease density (using the inverse move) by a bilipschitz map whose
  bilipschitz bound is bounded in terms of $D$.
  
\begin{figure}[ht]
    \centering
\includegraphics[width=.5\hsize]{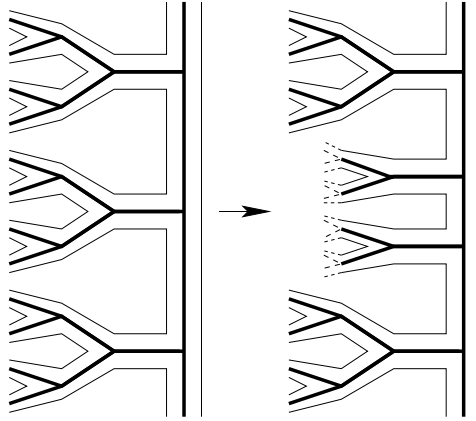}
    \caption{Increasing number of vertices along $\partial_0X$ by a 
    \emph{depth one splitting}}
    \label{Figure 1}
  \end{figure}
 
  One can then apply this process simultaneously on disjoint 
  intervals to change the local density along disjoint intervals. For 
  instance, applying  the above doubling procedure to all vertices 
  along $\partial_0X$ doubles the
  density, further, since it affects disjoint bounded sets, we still 
  have bilipshitz bound $k$. Similarly, given two disjoint intervals, 
  one could, for instance, 
  increase the local density by a factor of $D$ on one of the 
  intervals and decrease it on the other by a different factor $D'$; 
  since the intervals are disjoint the  
  bilipschitz bound depending only on the largest factor, which is 
  uniformly bounded by $K'$, the bilipschitz bound for $\Phi_{0}$.

  By these means we can, by replacing $X$ by its image under a
  bilipshitz map with bilipshitz bound bounded in terms of $K'$ and
  which is an isometry on $\partial_0X$, assure that the number of
  vertices along $\partial_0X$ and $\partial_0X_0$ matches to within a
  fixed constant over any interval in $\partial_0X$ and the
  corresponding image under $\Phi_0$. We now construct a bilipshitz
  map from this new $X$ to $X_0$ by first extending $\Phi_0$ to a
  $1$-neighborhood of $\partial_0X$ and then extending over the rest
  of $X$ by isometries of the components of the complement of this
  neighborhood.  By composing the two bilipshitz maps we get a 
  bilipshitz homeomorphism $\Psi$ from the original $X$ that does what
  we want on $\partial_0X$ while on every other boundary component
  $\partial_iX$ it is an isometry outside an interval of length
  bounded in terms of $K'$.

  Now choose arbitrary $K>1$. On $\partial_iX$ we can find an
  interval $J$ of length bounded in terms of $K'$ and $K$ that
  includes the interval $J_0$ on which our map is not an isometry, and
  whose length increases or decreases under $\Psi$ by a factor of at
  most $K$ (specifically, if the length of $J_0$ was multiplied by
  $s$, choose $J$ of length $\lambda\ell(J_0)$ with
  $\lambda\ge\max(\frac{K-Ks}{K-1},\frac{s-1}{K-1})$). Let $\Psi'$ be
  the map of $J$ that is a uniform stretch or shrink by the same
  factor (so the images of $\Psi'$ and $\Psi|J$ are identical). The
  following self-map $\alpha$ of a collar neighborhood $J\times[0,\epsilon]$
  restricts to $\Psi'\circ\Psi^{-1}$ on the left boundary $J\times \{0\}$
  and to the identity on the rest of the boundary:
  $$\alpha(x,t)=\frac{\epsilon-t}\epsilon\Psi'\circ\Psi^{-1}(x)+\frac
  t\epsilon x.$$ This $\alpha$ has bilipshitz constant bounded in terms
  of the bound on the left boundary and the length of $J$, hence
  bounded in terms of $K$ and $K'$. By composing $\Psi$ with $\alpha$
  on a collar along the given interval we adjust $\Psi|\partial_iX$
  to be a uniform stretch or shrink along this interval. We can do
  this on each boundary component other than $\partial_0X_0$. The
  result is a bilipschitz homeomorphism whose bilipschitz bound $L$ is
  still bounded in terms of $K'$ and $K$ and which satisfies the
  conditions of the theorem.
\end{proof}

We now deduce an analogue of the above Theorem in the case where the
boundary curves $\partial X$ are each labelled by one of a finite
number of colors, $C$, and the maps are required to be color
preserving. We call a labelling a \emph{bounded coloring} if there is
a uniform bound, so that given any point in $X$ and any color there is
a boundary component of that color a uniformly bounded distance away.
The lift of a coloring on a compact surface yields a bounded coloring.
We now fix a bounded coloring on our ``standard model'' $X_{0}$,
further, we choose this coloring so that it satisfies the following
regularity condition which is stronger than the above hypothesis: for
every point on a boundary component and every color in $C$, there is
an adjacent boundary component with that color a bounded distance from
the given point.  Call the relevant bound $B$.

\begin{theorem}\label{coloredfattrees}
  Let $X$ be as above with a chosen boundary component $\partial_0X$ 
  and fix a bounded coloring on the elements of $\partial X$. 
  Then there exists $K$ and a function $\phi\colon \R\to \R$ such that
  for any $K'$ and any color preserving 
  $K'$-bilipschitz homeomorphism $\Phi_0$ from
  $\partial_0X$ to a boundary component $\partial_0X_0$ of the
  ``standard model'' $X_0$, $\Phi_0$ extends to a
  $\phi(K')$-bilipschitz homeomorphism $\Phi\colon X\to X_0$ 
   which is $K$-bilipschitz on every other boundary component and 
   which is a color preserving map from $\partial X$ to $\partial 
   X_{0}$.
\end{theorem}

\begin{proof}
  As in the proof of Theorem~\ref{fattrees} we may assume $X$ is
  isometric to our standard model $X_{0}$ and then we proceed in two
  steps, first extending near $\partial_{0} X$, then extending over
  the rest of $X$.
    
  To extend near $\partial_{0} X$ we proceed as in the proof of
  Theorem~\ref{fattrees}, except now we need to not only match
  density, but colors as well. Instead of using only a \emph{depth one
    splitting} as in Figure~\ref{Figure 1}, one may perform a
  \emph{depth~$n$ splitting} by choosing a vertex at distance~$n$ from
  $\partial_{0} X$ and then moving that vertex to be adjacent to
  $\partial_{0} X$; this bilipschitz map increases the density of
  vertices along a given boundary component. Note that a depth~$n$
  move (and its inverse) can be obtained as a succession of depth~1
  moves and their inverses, so using such moves is only to yield a
  more concise language. Since the coloring of $X$ is a bounded
  coloring, from any point on $\partial_{0} X$, there is a uniform
  bound on the distance to a vertex adjacent to a boundary component
  of any given color. Thus, with a bounded bilipschitz constant we may
  alter the density and coloring as needed.
    
  As in the previous proof, we may extend to a map which does what is
  required on $\partial_0X$ and which is an isometry, but not
  preserving boundary colors, outside a neighborhood of $\partial_0X$,
  and which is a $K$-bilipschitz map on the boundary components other
  than $\partial_0X$ with $K$ close to $1$.

  Step two will be to apply a further bilipschitz map that fixes up
  colors on these remaining boundary components.

\begin{figure}[ht]
    \centering
    \labellist\small\hair 2.5pt
    \pinlabel{$\partial_{1} X$} at -10 0
    \pinlabel{$\partial_{0} X$} at 30 25
    \pinlabel{$\partial_{2} X$} at 100 25
    \pinlabel{$\partial_{3} X$} at 170 25
    \endlabellist
\includegraphics[width=.5\hsize]{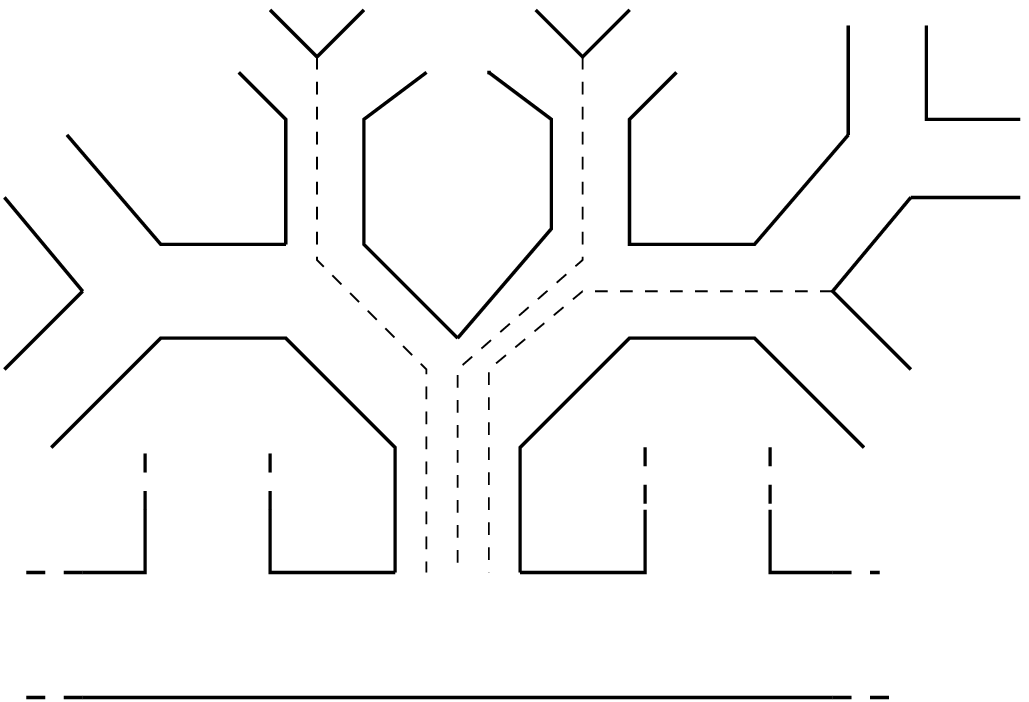}
\vglue10pt
\includegraphics[width=.5\hsize]{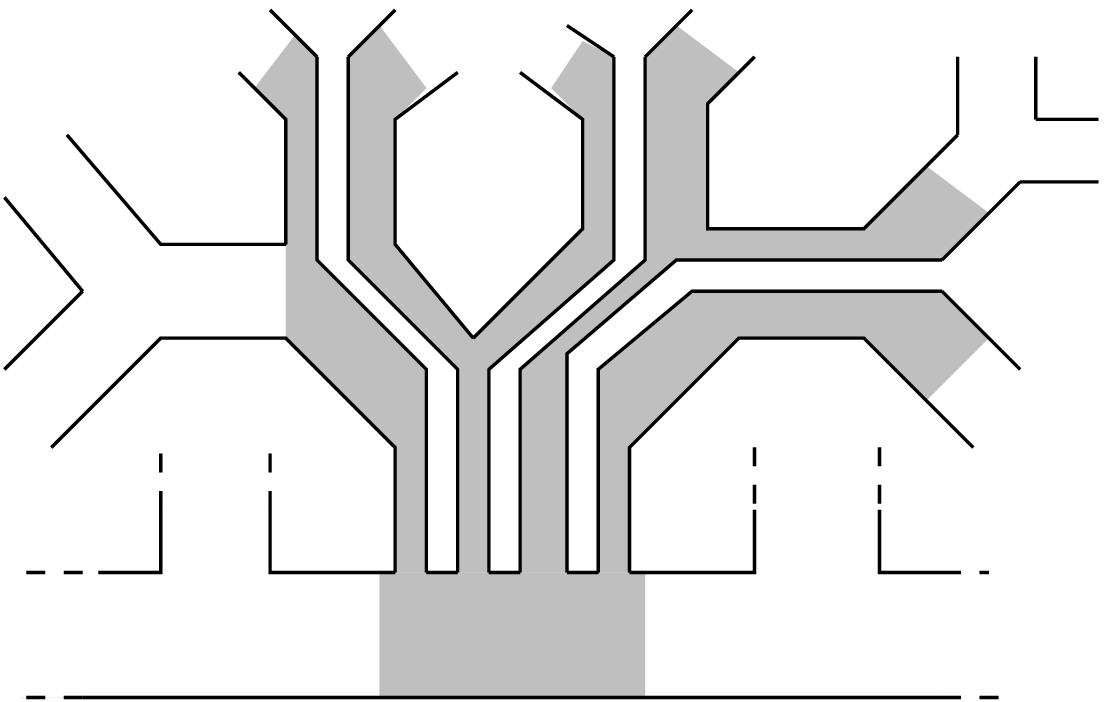}
\vglue10pt
\includegraphics[width=.75\hsize]{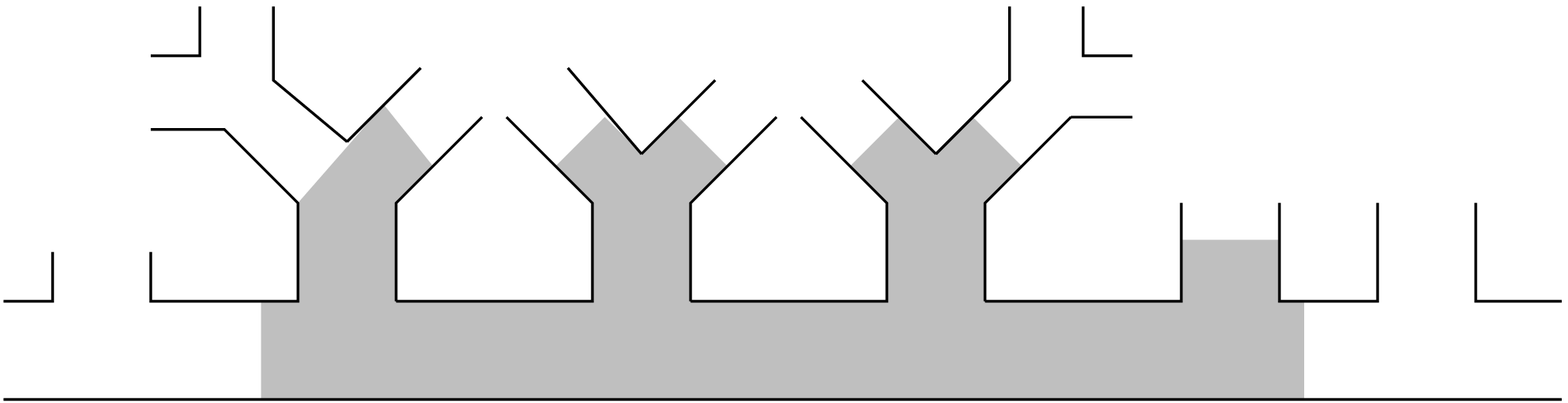}
    \caption{Adjusting colors along $\partial_1X$. Three depth 3 moves
    are illustrated. The shaded region shows where metric is
    adjusted. The final metric is shown at the bottom. In the first
    and last pictures all edges should have the same length since both
    are isometric to the standard model (some
    distortion was needed to draw them).}
    \label{Figure 2}
  \end{figure}

  Consider a boundary component $\partial_1X$ adjacent to
  $\partial_0X$. We want to make colors correct on boundary components
  adjacent to $\partial_1X$. They are already correct on $\partial_0X$
  and the boundary components adjacent on each side of this.  Call
  these $\partial_2X$ and $\partial'_2X$. As we move along
  $\partial_1X$ looking at boundary components, number the boundary
  components $\partial_0X$, $\partial_2X$, $\partial_3X,\dots$ until
  we come to a $\partial_{j+1}X$ which is the wrong color.  We will
  use splitting moves to bring new boundary components of the desired
  colors in to be adjacent to $\partial_1X$ between $\partial_jX$ and
  $\partial_{j+1}X$. By our regularity assumption on $X_0$ we will
  need to add at most $2B$ new boundary components before the color of
  $\partial_{j+1}X$ is needed; thus we need to perform at most $2B$
  splitting moves.  Moreover, the bounded coloring hypothesis implies
  that each of these splitting moves can be chosen to be of a
  uniformly bounded depth (note that the bounded coloring assumption
  implies that at any point of $X$ and any direction in the underlying
  tree there will be any desired color a uniformly bounded distance
  away).  We repeat this process along all of $\partial_1X$ in both
  directions to make colors correct. The fact that we do at most $2B$
  such moves for each step along $\partial_1X$ means that we affect
  the bilipschitz constant along $\partial_1X$ by at most a factor of
  $2B+1$. Since bounded depth splitting moves have compact support and
  since there are at most $2B$ of these performed between any pair
  $\partial_{j}X$, $\partial_{j+1}X$, we see that the bilipschitz
  constant need to fix this part of $\partial_{1} X$ is bounded in
  terms of $B$ and the bounded coloring constant. Since for $i\neq j$
  the neighborhoods affected by fixing the part of $\partial_{1}X$
  between $\partial_{j}X$, $\partial_{j+1}X$ are disjoint from those
  affected by fixing between $\partial_{i}X$, $\partial_{i+1}X$, we
  see that fixing all of $\partial_{1}X$ requires a bilipschitz bound
  depending only on $B$ and the bounded coloring constant; let us call
  this bilipschitz bound $C$.

  We claim that repeating this process boundary component by
  boundary component one can keep the bilipschitz constant under control and 
  thus prove the theorem. To see this, consider Figure \ref{Figure 2},
  which illustrates a typical set of splitting moves and shows the
  neighborhoods on which the metric has been altered. Since we can
  assure an upper
    bound on the diameter of these neighborhoods, we can
  bound the bilipshitz constants of the modifications, this is the 
  constant $C$ above. Only a bounded
  number of the neighborhoods needed for later modifications will
  intersect these neighborhoods, this yields  an overall bilipschitz 
  bound  which is at most a bounded  power of $C$.
\end{proof}

\section{Closed graph manifolds}
The following theorem answers Question~1.2 of
Kapovich--Leeb \cite{KapovichLeeb:3manifolds}. It is a special case of
Theorem \ref{th:main}, but we treat it here separately since its proof
is simple and serves as preparation for the general result.
\begin{theorem}\label{graphsqi} Any two closed non-geometric graph
  manifolds have bilipschitz homeomorphic universal covers. In
  particular, their fundamental groups are quasi-isometric.
\end{theorem}

Let us begin by recalling:
\begin{lemma}[Kapovich-Leeb \cite{KapovichLeeb:3manifolds};
    Neumann \cite{Neumann:commens}] 
    \label{kapovichleeb:lemma1}
    Any non-geometric graph manifold has an orientable finite cover where all
    Seifert components are circle bundles over orientable
    surfaces of genus~$\geq 2$.  Furthermore, one can arrange that the
    intersection numbers of the fibers of adjacent Seifert components
    are $\pm 1$.
\end{lemma}
If we replace our graph manifold by a finite cover as in the above
lemma then we have a trivialization of the circle bundle on the
boundary of each Seifert piece using the section given by a fiber of a
neighboring piece. The fibration of this piece then has a relative
Euler number.
\begin{lemma}[Kapovich--Leeb \cite{KapovichLeeb:3manifolds}]
\label{kapovichleeb:lemma2}
Up to a bilipschitz homeomorphism of the universal cover, we can
assume all the above relative Euler numbers are $0$.
\end{lemma}
A graph manifold $G$ as in the last lemma is what Kapovich and Leeb
call a ``flip-manifold.'' It is obtained by gluing together finitely
many manifolds of the form (surface)$\times S^1$ by gluing using maps
of the boundary tori that exchange base and fiber coordinates. We can
give it a metric in which every fiber $S^1$ (and hence every boundary
circle of a base surface) has length $1$.

A topological model for the universal cover $\tilde G$ can be obtained by
gluing together infinitely many copies of $X_0\times \R$ according to
a tree, gluing by the ``flip map'' \begin{tiny}$\begin{pmatrix}
  0&1\\1&0
\end{pmatrix}$\end{tiny}$\colon\R\times\R\to \R\times\R$ when gluing
boundary components. We call the resulting manifold $Y$. 

We wish to show that $\tilde G$ is bilipschitz homeomorphic to $Y$.

\begin{proof}[Proof of Theorem~\ref{graphsqi}]
  The universal cover of each Seifert component of $G$ is identified
  with $\tilde{S_i}\times \R$, where $S_i$ is one of a finite
  collection of compact surfaces with boundary. Choose a number $K$
  sufficiently large that Theorem \ref{fattrees} applies for each of
  them. Choose a bilipschitz homeomorphism from one piece $\tilde 
  S_i\times \R$
  of $\tilde G$ to a piece $X_0\times \R$ of $Y$, preserving the
  (surface)$\times \R$ product structure. We want to extend to a
  neighboring piece of $\tilde G$. On the common boundary $\R\times
  \R$ we have a map that is of the form $\phi_1\times \phi_2$ with
  $\phi_1$ and $\phi_2$ both bilipschitz. By Theorem \ref{fattrees} we can
  extend over the neighboring piece by a product map, and on the other
  boundaries of this piece we then have maps of the form
  $\phi'_1\times \phi_2$ with $\phi'_1$ $K$-bilipschitz. We do this
  for all neighboring pieces of our starting piece. Because of the
  flip, when we extend over the next layer we have maps on the outer
  boundaries that are $K$-bilipschitz in both base and fiber. We can
  thus continue extending outwards inductively to construct our
  desired bilipschitz map.
\end{proof}

\section{Graph manifolds with boundary} \label{S:Proof}

A non-geometric graph manifold $M$ has a minimal decomposition along
tori and Klein bottles into geometric (Seifert fibered) pieces, called
the \emph{geometric decomposition}. The cutting surfaces are then
$\pi_1$ injective. In this decomposition one cuts along one-sided
Klein bottles; this differs from JSJ, where one would cut along the
torus neighborhood boundaries of these Klein bottles.  (See, e.g.,
Neumann--Swarup \cite{neumann-swarup} Section 4.)

We associate to this decomposition its \emph{decomposition graph},
which is the graph with a vertex for each Seifert component of $M$ and
an edge for each decomposing torus or Klein bottle. If there are no
one-sided Klein bottles then this graph is the graph of the associated
graph of groups decomposition of $\pi_1(M)$.  (If there are
decomposing Klein bottles, the graph of groups has, for each Klein
bottle, an edge to a new vertex rather than a loop. This edge
corresponds to an amalgamation to a Klein bottle group along a
$\Z\times\Z$, and corresponds also to an inversion for the action of
$\pi_1(M)$ on the Bass-Serre tree. Using a loop rather than an edge
makes the Bass-Serre tree a weak covering of the decomposition graph.)

We color vertices of the decomposition graph \underline{\bf b}lack or
\underline{\bf w}hite according to whether the Seifert piece includes
a boundary component of $M$ or not (\underline{\bf b}ounded or
\underline{\bf w}ithout boundary).

A second graph we consider is the  \bicolor{}ed decomposition graph
for the decomposition of the universal cover $\tilde M$ into its
fibered pieces. We denote it $BS(M)$ and call it the \emph{\bicolor{}ed
  Bass-Serre tree,} since it is the Bass-Serre tree for our graph of
groups decomposition.  It can be obtained from the \bicolor{}ed
decomposition graph by replacing each edge by a countable infinity of
edges between its endpoints, and then taking the universal cover of
the resulting graph.

A \emph{weak covering map} from a \bicolor{}ed graph $\Gamma$ to a
\bicolor{}ed graph $\Gamma'$ is a color-preserving graph homomorphism
$\phi\colon\Gamma\to \Gamma'$ with the property that for any vertex
$v$ of $\Gamma$ and every edge $e'$ at $\phi(v)$, there is at least
one edge $e$ at $v$ mapping to $e'$.  An example of such a map is
the map that collapses any multiple edge of $\Gamma$ to a single edge.
Any covering map of non-geometric graph
manifolds induces a weak covering map of their \bicolor{}ed decomposition
graphs. 

Note that if a weak covering map exists from $\Gamma$ to $\Gamma'$
then $\Gamma$ and $\Gamma'$ will have isomorphic \bicolor{}ed Bass-Serre
trees.  The equivalence relation on \bicolor{}ed graphs generated by the
relation of existence of a weak covering map will be called
\emph{bisimilarity}. We shall prove in the next section:

\begin{proposition}\label{prop:minimal and Bass Serre}
  If we restrict to countable connected graphs then each equivalence
  class of \bicolor{}ed graphs includes two characteristic elements: a
  unique tree that weakly covers every element in the class (the
  Bass-Serre tree); and a unique minimal element, which is weakly
  covered by all elements in the class.
\end{proposition}
For example, if all the vertices of a graph have the same color, then
the minimal graph for its bisimilarity class is a single vertex with a
loop attached and the Bass-Serre tree is the single-colored regular
tree of countably infinite degree.

Our main theorem is:
\begin{theorem}\label{th:main}  
  If $M$ and $M'$ are non-geometric graph manifolds then the
  following are equivalent:
  \begin{enumerate}
  \item\label{it:m1} $\tilde M$ and $\tilde M'$ are bilipschitz
    homeomorphic.
  \item\label{it:m2} $\pi_1(M)$ and $\pi_1(M')$ are
    quasi-isometric.
  \item\label{it:m3} $BS(M)$ and $BS(M')$ are isomorphic as
    \bicolor{}ed trees.
  \item\label{it:m4} The minimal \bicolor{}ed graphs in the bisimilarity
    classes of the decomposition graphs $\Gamma(M)$ and $\Gamma(M')$
    are isomorphic.
  \end{enumerate}
\end{theorem}
\begin{proof}
Clearly (\ref{it:m1}) implies (\ref{it:m2}). The equivalence of
(\ref{it:m3}) and (\ref{it:m4}) is Proposition \ref{prop:minimal and
  Bass Serre}.
Kapovich and Leeb~\cite{KapovichLeeb:3manifolds} proved that any
quasi-isometry essentially preserves the geometric decomposition of
Haken manifolds, and therefore induces an isomorphism between their
Bass-Serre trees. To prove the theorem it remains to show that
(\ref{it:m3}) or (\ref{it:m4}) implies (\ref{it:m1}).

Suppose therefore that $M$ and $M'$ are non-geometric graph manifolds
that satisfy the equivalent conditions (\ref{it:m3}) and
(\ref{it:m4}). Let $\Gamma$ be the minimal graph in the bisimilarity
class of $\Gamma(M)$ and $\Gamma(M')$. It suffices to show that each
of $\tilde M$ and $\tilde M'$ is bilipschitz homeomorphic to the
universal cover of some standard graph manifold associated to
$\Gamma$. There is therefore no loss in assuming that $M'$ \emph{is}
such a standard graph manifold; ``standard'' will mean that
$\Gamma(M')=\Gamma$ and that each loop at a vertex in $\Gamma$
corresponds to a decomposing Klein bottle (i.e., a boundary torus of
the corresponding Seifert fibered piece that is glued to itself by a
covering map to the Klein bottle).

Denote the set of pairs consisting of a vertex of $\Gamma$ and an
outgoing edge at that vertex by $C$. Since the decomposition graphs
$\Gamma(M)$, $\Gamma(M')$, $BS(M)$, and $BS(M')$ for $M$, $M'$,
$\tilde M$, and $\tilde M'$ map to $\Gamma$, we can label the boundary
components of the geometric pieces of these manifolds by elements of $C$.

Our desired bilipschitz map can now be constructed inductively as in 
the proof of Theorem~\ref{graphsqi}, at each stage of the process 
having extended over some submanifold $Y$ of $\tilde G$. The 
difference from the situation there is that now when we extend the map 
from $Y$ over a further fibered piece $X\times \R$, we must
make sure that we are mapping boundary components to boundary
components with the same $C$--label. That this can be done is exactly 
the statement of Theorem~\ref{coloredfattrees}.
\end{proof}

\begin{remark} With some work, Theorem~\ref{th:main} can be generalized to 
cover many situations outside of the context of $3$-manifolds; such a 
formulation will appear in a forthcoming paper.    
\end{remark}

\section{ \Bicolor{}ed graphs}\label{sec:bicolor}

\def\b{\bf b}\def\w{\bf w}
\begin{definition}\label{def:bicolor}
  A \emph{graph} $\Gamma$ consists of a vertex set $V(\Gamma)$ and an
  edge set $E(\Gamma)$ with a map $\epsilon\co E(\Gamma)\to
  V(\Gamma)^2/C_2$ to the set of unordered pairs of
  elements of $V(\Gamma)$.

 A \emph{\bicolor{}ed graph} is a graph $\Gamma$ with a ``coloring'' $c\co
  V(\Gamma) \to \{\b,\w\}$. 

  A \emph{weak covering} of \bicolor{}ed graphs is a graph homomorphism
  $f\co \Gamma \to \Gamma' $ which respects colors and has the
  property: for each $v\in V(\Gamma)$ and for each edge $e'\in
  E(\Gamma')$ at $f(v)$ there exists an $e\in E(\Gamma)$ at $v$ with
  $f(e)=e'$.
\end{definition}
From now on, all graphs we consider will be assumed to be
connected. It is easy to see that a weak covering is then
surjective. The graph-theoretic results are valid for $n$--color
graphs, but we only care about $n=2$.
\begin{definition}\label{def:bisimilar}
  \Bicolor{}ed graphs $\Gamma_1,\Gamma_2$ are
  \emph{bisimilar}, written $\Gamma_1\sim\Gamma_2$, if $\Gamma_1$ and
  $\Gamma_2$ weakly cover some common \bicolor{}ed graph.
\end{definition}
The following proposition implies, among other things, that this
definition agrees with our earlier version.
\begin{proposition}\label{prop:bisimeq}
  The bisimilarity relation $\sim$ is an equivalence relation.
  Moreover, each equivalence class has a unique minimal element up to
  isomorphism.
\end{proposition}
\begin{lemma}\label{le:bisimeq}
  If a \bicolor{}ed graph $\Gamma$ weakly covers each of a collection of graphs
$\{\Gamma_i\}$ then the $\Gamma_i$ all weakly cover
some common $\Gamma'$.
\end{lemma}
\begin{proof}
  The graph homomorphism that restricts to a bijection on the vertex set but
  identifies multiple edges with the same ends to a single edge is a
  weak covering. Moreover, if we do this to both graphs $\Gamma$ and
  $\Gamma_i$ of a weak covering $\Gamma\to\Gamma_i$ we still have a
  weak covering. So there is no loss in assuming all our graphs have no
  multiple edges. A graph homomorphism $\Gamma\to\Gamma_i$ is then
  determined by its action on vertices. The induced equivalence
  relation $\equiv$ on vertices of $\Gamma$ satisfies the property:
\begin{itemize}
\item[] 
If $v\equiv v_1$ and $e$ is an
  edge with $\epsilon(e)=\{v,v'\}$ then there exists an edge $e_1$ with
$\epsilon(e_1)=\{v_1,v'_1\}$ and $v'\equiv v'_1$. 
\end{itemize}
Conversely, an equivalence relation on vertices of $\Gamma$ with this
property induces a weak covering. We must thus just show that if we
have several equivalence relations on $V(\Gamma)$ with this property,
then the equivalence relation $\equiv$ that they generate still has
this property. Suppose $v\equiv w$ for the generated relation. Then we
have $v=v_0\equiv_1 v_1\equiv_2 \dots \equiv_n v_n=w$ for some $n$,
where the equivalence relations $\equiv_i$ are chosen from our given
relations. Let $e_0$ be an edge at $v=v_0$ with other end at $v'_0$.
Then the above property guarantees inductively that we can find an
edge $e_i$ at $v_i$ for $i=1,2,\dots,n$, with other end at $v'_i$ and with
$v'_{i-1}\equiv_i v'_i$. Thus we find an edge $e_n$ at $w=v_n$ whose
other end $v'_n$ satisfies $v'_0\equiv v'_n$.
\end{proof}
\begin{proof}[Proof of Proposition \ref{prop:bisimeq}]
  We must show that $\Gamma_1\sim\Gamma_2\sim\Gamma_3$ implies
  $\Gamma_1\sim\Gamma_3$. Now $\Gamma_1$ and $\Gamma_2$ weakly cover a
  common $\Gamma_{12}$ and $\Gamma_2$ and $\Gamma_3$ weakly cover some
  $\Gamma_{23}$. The lemma applied to $\Gamma_2, \{\Gamma_{12},
  \Gamma_{23}\}$ gives a graph weakly covered by all three of
  $\Gamma_1,\Gamma_2,\Gamma_3$, so $\Gamma_1\sim\Gamma_3$. 

  The minimal element in a bisimilarity class is found by applying the
  lemma to an element $\Gamma$ and
  the set $\{\Gamma_i\}$ of all \bicolor{}ed graphs that $\Gamma$ weakly covers. 
\end{proof}

\begin{proposition} \label{prop:bicolor tree}
If we restrict to \bicolor{}ed graphs all of whose
  vertices have countable valence (so the graphs are also countable,
  by our connectivity assumption), then each bisimilarity class
  contains a tree $T$, unique up to isomorphism, that weakly covers
  every element of the class.  It can be constructed as follows: If\/
  $\Gamma$ is in the bisimilarity class, duplicate every edge of\/
  $\Gamma$ a countable infinity of times, and then take the universal
  cover of the result (in the topological sense).
\end{proposition}
Note that uniqueness of $T$ in the above proposition depends on the
fact that $T$ is a tree; there are many different \bicolor{}ed graphs 
that weakly cover every \bicolor{}ed graph in a given bisimilarity class.
\begin{proof}[Proof of Proposition \ref{prop:bicolor tree}]
  Given a \bicolor{}ed graph $\Gamma$, we can construct a
  tree $T$ as follows: Start with one vertex $x$, labeled by a vertex
  $v$ of $\Gamma$. Then for each vertex $w$ of $\Gamma$ connected to
  $v$ by an edge, add infinitely many edges at $x$ leading to vertices
  labeled $w$. Then repeat the process at these new vertices and
  continue inductively. Finally forget the $\Gamma$--labels on the
  resulting tree and only retain the corresponding
  $\{\b,\w\}$--labels.

  If $\Gamma$ weakly covers a graph $\Gamma'$, then
  using $\Gamma'$ instead of $\Gamma$ to construct the above tree $T$
  makes no difference to the inductive construction.  Thus $T$ is an
  invariant for bisimilarity. It clearly weakly covers the
  original $\Gamma$, and since $\Gamma$ was arbitrary in the
  bisimilarity class, we see that $T$ weakly covers anything in the
  class.

  To see uniqueness, suppose $T'$ is another tree that weakly covers
  every element of the bisimilarity class. Then $T'$ weakly covers the
  $T$ constructed above from $\Gamma$. Composing with $T\to \Gamma$
  gives a weak covering $f\co T'\to \Gamma$ for which infinitely many
  edges at any vertex $v\in V(T')$ lie over each edge at the vertex
  $f(v)\in V(\Gamma)$. It follows that $T'$ itself can be constructed
  from $\Gamma$ as in the first paragraph of this proof, so $T'$ is
  isomorphic to $T$.
\end{proof}
Using a computer we have found (in about 5 months of processor time):
\begin{proposition}
  The number of connected minimal \bicolor{}ed graphs with $n$
  vertices of which exactly $b$ are black (excluding the two
  $1$--vertex graphs with no edges) is given by the table:
 
\smallskip\centerline{\begin{tiny}\begin{tabular}{|l|rccccccc|c|}
\hline
$n$ &$b:0$&$1$&$2$&$3$&$4$&$5$&$6$&$7$&total\\
\hline
$1$&$1$&$1$&$0$&$0$&$0$&$0$&$0$&$0$&$2$\\
$2$&$0$&$4$&$0$&$0$&$0$&$0$&$0$&$0$&$4$\\
$3$&$0$&$10$&$10$&$0$&$0$&$0$&$0$&$0$&$20$\\
$4$&$0$&$56$&$61$&$56$&$0$&$0$&$0$&$0$&$173$\\
$5$&$0$&$446$&$860$&$860$&$446$&$0$&$0$&$0$&$2612$\\
$6$&$0$&$6140$&$17084$&$20452$&$17084$&$6140$&$0$&$0$&$66900$\\
$7$&$0$&$146698$&$ 523416$&$ 755656$&$ 755656$&$ 523416$&$
146698$&$0$&$2851540$\\
$8$&$0$&$6007664$&$25878921$&$44839104$&$48162497$&$44839104$&$25878921$&$6007664$&$201613875$\\
\hline
\end{tabular}\end{tiny}}\par%\qed    
\end{proposition}
The proposition shows, for example, that there are 199 quasi-isometry
classes for non-geometric graph manifolds having four or fewer Seifert
pieces ($199=2+4+20+173$). In the next subsection we list the
corresponding 199 graphs. These were found by hand before programming
the above count. This gives some confidence that the computer program is
correct.

\subsection{Enumeration of minimal \bicolor{}ed graphs up to 4 vertices}
We only consider connected graphs and we omit the two 1-vertex graphs with no
edges. In the following table ``number of graphs $n+n$'' means $n$
graphs as drawn and $n$ with $\b$ and $\w$ exchanged. Dotted loops in
the pictures represent loops that may
or may not be present and sometimes carry labels $x,x',\dots$
 referring to the two-element set
$\{$``present'', 
``absent''$\}$.

\def\lloop#1{\ar@(ul,dl)@{#1}[]}
\def\rloop#1{\ar@(ur,dr)@{#1}[]}
\def\dloop#1{\ar@(dl,dr)@{#1}[]}
\def\uloop#1{\ar@(ul,ur)@{#1}[]}
$$\xymatrix@R=10pt@C=36pt@M=0pt@W=0pt@H=0pt{
&&&&\hbox{number of graphs}\\
\hbox{1 vertex:}
&\Dot\rloop-&\Circ\rloop-&&2\\\\\\
\hbox{2 vertices:}
&\Dot\lineto[r]\lloop.&\Circ\rloop.&&4\\\\\\
\hbox{3 vertices:}
&\Dot\lloop.\lineto[r]&\Circ\dloop.\lineto[r]&\Circ\rloop.&8+8\\\\
&\Circ\lineto[r]&\Dot\dloop.\lineto[r]&\Circ\rloop-&2+2
\hbox to 0 pt{\qquad (total: 20)\hss}\\\\\\
}$$
$$\xymatrix@R=10pt@C=36pt@M=0pt@W=0pt@H=0pt{
\hbox to 0pt{\hss 4 vertices:\quad}
&\Dot\lloop.\lineto[r]&\Circ\dloop.\lineto[r]&\Circ\dloop.\lineto[r]&\Circ\rloop.&16+16\\\\
&\Circ\lloop.\lineto[r]&\Dot\dloop.\lineto[r]&\Circ\dloop.\lineto[r]&\Circ\rloop.&16+16\\\\
&\Dot\lloop.\lineto[r]&\Dot\dloop.\lineto[r]&\Circ\dloop.\lineto[r]&\Circ\rloop.&16\\\\
%&\Dot\lloop.\lineto[r]&\Dot\dloop.\lineto[r]&\Dot\dloop.\lineto[r]&\Circ\rloop.&16\\\\
&\Dot\lineto[r]&\Circ\dloop.\lineto[r]&\Circ\dloop.\lineto[r]&\Dot\rloop-&4+4\\\\
&\Dot\lloop._(.75)x\lineto[r]&\Circ\dloop._(.75)y\lineto[r]&\Dot\dloop._(.75){x'}\lineto[r]&\Circ\rloop.^(.75){y'}&12\hbox
to 12pt{\quad($x\ne x'$ or $y\ne y'$)\hss} \\\\
%}$$ $$\xymatrix@R=10pt@C=36pt@M=0pt@W=0pt@H=0pt{
&&&&\Circ\\
&&\Dot\lloop.\lineto[r]&\Dot\uloop.\lineto[ur]\lineto[dr]&&4+4\\
&&&&\Circ\rloop-\\\\
&&&&\Circ\\
&&\Dot\lloop.\lineto[r]&\Circ\uloop.\lineto[ur]\lineto[dr]&&4+4\\
&&&&\Circ\rloop-\\\\
&&&&\Circ\rloop.\lineto[dd]\\
&&\Circ\lloop.\lineto[r]&\Circ\uloop.\lineto[ur]\lineto[dr]&&16+16\\
&&&&\Dot\rloop.\\\\
}$$ $$\xymatrix@R=10pt@C=36pt@M=0pt@W=0pt@H=0pt{
&&&&\Circ\rloop.\lineto[dd]\\
&&\Dot\lloop._(.75)x\lineto[r]&\Circ\uloop.\lineto[ur]\lineto[dr]&&8+8\hbox
to 0 pt{\quad ($x\ne x'$)\hss}\\
&&&&\Dot\rloop.^(.75){x'}\\\\
%}$$ $$\xymatrix@R=10pt@C=36pt@M=0pt@W=0pt@H=0pt{
%&&&&\Circ\rloop.\lineto[dd]\\
%&&\Dot\lloop-\lineto[r]&\Circ\uloop.\lineto[ur]\lineto[dr]&&4+4\\
%&&&&\Dot\\\\
&&\Dot\lloop.\lineto[r]\lineto[dd]&\Circ\rloop-\lineto[dd]\\
&&&&&4+4\\
&&\Circ\lineto[r]&\Circ\rloop.\\\\
&&\Dot\lineto[r]\lineto[dd]&\Circ\rloop-\lineto[dd]\\
&&&&&1\\
&&\Circ\lineto[r]&\Dot\rloop-\\\\
&&&&&\hbox to 36 pt{Total for 4 vertices: 173\hss}\\\\
}$$
\subsection{Algorithm for finding the minimal \bicolor{}ed graph}

\newcommand{\Adj}{\operatorname{Adjacent}}
\newcommand{\Cur}{\operatorname{CurrentColor}}
\newcommand{\MaxC}{\operatorname{MaxColor}}
Let $\Gamma$ be a connected \bicolor{}ed graph. We wish to construct the
minimal \bicolor{}ed graph $\Gamma_0$ for which there is a weak covering
$\Gamma\to\Gamma_0$. Note that any coloring $c\colon V(\Gamma)\to C$
of the vertices of $\Gamma$ induces a graph homomorphism to a graph
$\Gamma_c$ with vertex set $C$ and with an edge connecting the
vertices $w_1,w_2\in C$ if and only if there is some edge connecting a
$v_1,v_2\in V(\Gamma)$ with $c(v_i)=w_i$, $i=1,2$.

We start with $C$ containing just our original two colors, which we
now call $0,1$, and gradually enlarge $C$ while modifying $c$ until
the the map $\Gamma\to \Gamma_c$ is a weak covering. For a vertex $v$
let $\Adj(v)$ be the set of colors of vertices connected to $v$ by an
edge (these may include $v$ itself). We shall always call our coloring
$c$, even as we modify it.

\begin{enumerate}
\item $\Cur=0$; $\MaxC=1$;
\item {\bf While} $\Cur\le\MaxC$;
  \begin{enumerate}
  \item {\bf If} there are two vertices $v_1,v_2$ with $c(v_i)=\Cur$
    that have different $\Adj(v_i)$'s;
  \item {\bf Then} increment $\MaxC$ and add it to the set $C$, change the
    color of each $v$ with $c(v)=\Cur$ and $\Adj(v)=\Adj(v_1)$ to
    $\MaxC$, and then set $\Cur=0$;
  \item {\bf Else} increment $\Cur$;
  \item {\bf End If};
  \end{enumerate}
\item {\bf End While}.
\end{enumerate}

We leave it to the reader to verify that this algorithm terminates
with $\Gamma\to\Gamma_c$ the weak covering to the minimal \bicolor{}ed
graph (in step (2b) we could add a new color for each new value of
$\Adj(v)$ with $v\in\{v:c(v)=\Cur\}$ rather than for just one of them;
this seems \emph{a priori} more efficient but proved hard to program
efficiently). The algorithm is inspired by Brendan McKay's ``nauty''
\cite{mckay}; we are grateful to Dylan Thurston for the suggestion.

Counting the number of minimal \bicolor{}ed graphs with $b$ black
vertices and $w$ white vertices is now easy. We order the vertices
$1,\dots,b,\dots,b+w$ and consider all connected graphs on this vertex
set. For each we check by the above procedure if it is minimal and if
so we count it.  Finally, we divide our total count by $b!w!$ since
each graph has been counted exactly that many times (a minimal
\bicolor{}ed graph has no automorphisms).
\section{Artin groups}

An \emph{Artin group} is a group given by a presentation of the
following form:
$$
A = \langle x_1,..., x_n\mid (x_i,x_j)_{m_{ij}} = (x_j,x_i)_{m_{ji}}
\rangle
$$
where, for all $i \ne j $ in $\{1,\ldots,n\}$, $m_{ij} = m_{ji} \in
\{2,3,\ldots, \infty\}$ with $(x_i,x_j)_{m_{ij}} =x_ix_jx_i...$
($m_{ij}$ letters) if $m_{ij}<\infty$ and when $m_{ij}=\infty$ 
we do not add a defining relation between $x_{i}$ and $x_{j}$.
A concise way to present such a group is as a finite graph labeled by
integers greater than~$1$: such a graph has $n$ vertices, one for each
generator, and a pair of vertices are connected by an edge labeled by
$m_{ij}$ if $m_{ij}<\infty$.

An important class of Artin groups is the class of
\emph{right-angled Artin group}. These are Artin groups with each 
$m_{ij}$ either $2$ or $\infty$, i.e., the only
defining relations are commutativity relations between pairs of 
generators. These groups interpolate between the free group 
on~$n$ generators ($n$ vertices and no edges) and 
$\Z^{n}$ (the complete graph on~$n$ vertices).

We shall call a presentation tree \emph{big} if it has diameter $\ge
3$ or has diameter $2$ and at least one weight on it is $>2$. An Artin group
given by a non-big tree has infinite center and is
virtually (free)$\times \Z$.
The Artin groups given by non-big presentation trees thus fall into three
quasi-isometry classes ($\Z$, $\Z^2$, $F_2\times \Z$, where
$F_2$ is the 2-generator free group) and are not quasi-isometric to
any Artin group with big presentation trees (this follows, for 
instance, from \cite{KapovichLeeb:haken}).
We shall therefore only be concerned with Artin groups whose presentation
trees are big. For right-angled Artin groups
this just says the presentation tree has diameter larger than~$2$.

We use the term \emph{tree group} to refer to any 
Artin group whose presentation graph is 
a big tree. Any right-angled tree 
group is the fundamental group of a flip 
graph manifold: this is seen by identifying each diameter~$2$ region 
with a $\mbox{(punctured surface)}\times {\mathbb S}^{1}$ and noting 
that pairs of  
such regions are glued together by switching fiber and base 
directions.

Since any right-angled tree group corresponds to a graph manifold with boundary 
components in each Seifert piece, Theorem~\ref{th:main} yields 
immediately the following answer to
Bestvina's question about their quasi-isometry classification:  
\begin{theorem}\label{trees qi}
    Any pair of right-angled tree groups are quasi-isometric. 
\end{theorem}

This raises the following natural question:
\begin{question}\label{question:treegroups} 
    When is a finitely generated group $G$ quasi-isometric to a 
    right-angled tree group?
\end{question}
The simple answer 
is that  $G$ must be weakly commensurable with the fundamental group of a
non-geometric graph manifold with boundary components in every Seifert
component, this follows from our Theorem~\ref{th:main} and Kapovich--Leeb's  
quasi-isometric rigidity result for non-geometric $3$-manifolds 
\cite{KapovichLeeb:haken}.
But it is natural to ask the question within the class of
Artin groups, where this answer is not immediately helpful. We give 
the following answer, which in particular shows that right-angled 
tree groups are quasi-isometrically rigid in the class of 
right-angled Artin groups.
\begin{theorem}\label{treegroupsclassification} Let $G'$ be any Artin
  group and let $G$ be a right-angled tree group.  Then $G'$ is
  quasi-isometric to $G$ if and only if $G'$ has presentation graph a
  big even-labeled tree with all interior edges labeled 2. (An
  ``interior edge'' is an edge that does not end in a leaf of the
  tree.)
\end{theorem}

We first recall two results relevant to Artin groups given by trees. The
first identifies which Artin groups are $3$-manifold groups and the
second says what those $3$-manifolds are. 

\begin{theorem}[Gordon; \cite{Gordon:Coherence}]
    \label{artin manifolds}
The following are equivalent for an Artin group $A$:
\begin{enumerate}
\item $A$ is virtually a $3$-manifold group.
\item $A$ is a $3$-manifold group.
\item Each connected 
    component of its presentation graph is either a tree or a 
    triangle with each edge labeled~$2$.
\end{enumerate}
\end{theorem}
\begin{theorem}[Brunner \cite{Brunner}, Hermiller-Meier
  \cite{HermillerMeier}] 
\label{graph links}
  The Artin group associated to a weighted tree $T$ is the fundamental
  group of the complement of the following connected sum of
  torus links. For each $n$-weighted edge of $T$ associate a copy of
  the $(2,n)$-torus link and if $n$ is even associate each end of the
  edge with one of the two components of this link; if $n$ is odd
  associate both ends of the edge with the single component (a
  $(2,n)$-knot). Now take the connected sum of all these links, doing
  connected sum whenever two edges meet at a vertex, using the
  associated link components to do the sum.
\end{theorem}

(In Theorem \ref{graph links} the fact that for an odd-weighted edge
the $(2.n)$ torus knot can be associated with either end of the edge
shows that one can modify the presentation tree without changing the
group. This is a geometric version of the ``diagram twisting'' of
Brady, McCammond, M\"uhlherr, Neumann \cite{Brady+:Rigidity}.)

\begin{proof}[Proof of Theorem~\ref{treegroupsclassification}]
  Let $G'$ be an Artin group that is quasi-isometric to a right-angled
  tree group.  
  Right-angled tree groups are one-ended and hence $G$, and thus $G'$ 
  as well, is not freely decomposable.
  Thus the presentation graph for $G'$ is connected.
  
  By the quasi-isometric rigidity Theorem for $3$-manifolds, as stated in 
  the introduction, we know that $G'$ is weakly commensurable to a 
  $3$-manifold group.
     
  Unfortunately it is not yet known if every Artin group is torsion
  free. If we knew $G'$ is torsion free then we could argue as
  follows. 
  First, since $G'$ is torsion free, it follows that 
  $G'$ is commensurable
  with a 3-manifold group. Thus by Theorem \ref{artin manifolds} it
  \emph{is} a 3-manifold group and is a tree group. By Theorem
  \ref{th:main} the corresponding graph manifold must
  have boundary components in every Seifert component. Using Theorem
  \ref{graph links} it is then easy to see that this gives precisely
  the class of trees of the theorem. We say more on this in Theorem
  \ref{artin to graphmanifold} below.

  Since we only know that the quotient of $G'$ by a finite group,
  rather than $G'$ itself, is commensurable with a $3$-manifold group
  we cannot use Gordon's result (Theorem \ref{artin manifolds})
  directly. But we will follow its proof. 
  
  Gordon rules out most Artin groups being fundamental groups of
  $3$-manifolds by proving that they contain finitely generated
  subgroups which are not finitely presented (i.e., they are not
  \emph{coherent}).  Since Scott \cite{Scott:Coherent} proved
  $3$-manifold groups are coherent, and since coherence is a
  commensurability invariant, such Artin groups are not $3$-manifold
  groups.  Since coherence is also a weak commensurability invariant, 
  this also rules out these Artin groups in our situation.
    
  The remaining Artin groups which Gordon treats with a separate
  argument are those that include triangles with labels $(2,3,5)$ or
  $(2,2,m)$. The argument given by Gordon for these cases also applies
  for weak commensurability. (A simpler argument than Gordon's in the
  $(2,2,m)$ case is that $A$ then contains both a $\Z^3$ subgroup and
  a non-abelian free subgroup, which easily rules out weak
  commensurability with a $3$-manifold group.)
\end{proof}

The above argument leads also to the following generalization of
Gordon's theorem.
\begin{theorem}\label{gordon extension}
  An Artin group $A$ is quasi-isometric to a $3$-manifold group if and
  only if it is a $3$-manifold group (and is hence as in Theorem
  \ref{artin manifolds}).
\end{theorem}
\begin{proof} Fix an Artin group $A$ which is quasi-isometric to a 
    $3$-manifold group. By Papasoglu-Whyte \cite{PapasogluWhyte:ends}, 
    the reducible case 
    reduces to the irreducible case, so we assume the Artin group has a 
    connected presentation graph.
    
    The quasi-isometric rigidity Theorem for $3$-manifolds implies that $A$ is 
    weakly commensurable (or in some cases even commensurable) with a 
    $3$-manifold group, so as in the previous proof an easy modification of 
    Gordon's argument applies. 
\end{proof}

We can, in fact, more generally describe the quasi-isometry class of
any tree group $A$ in terms of Theorem \ref{treegroupsclassification}.
That is, we can describe the \bicolor{}ed decomposition graph for the graph
manifold $G$ whose fundamental group is $A$.
\begin{theorem}\label{artin to graphmanifold}
  The colored decomposition graph is obtained from the presentation
  tree of the Artin group by the following sequence of moves:
  \begin{enumerate}
\item Color all existing vertices black.
  \item For each odd-weighted edge, collapse the edge, thus
    identifying the vertices at its ends, and
    add a new edge from this vertex to a new leaf which is colored white.
  \item Remove any $2$-weighted edge leading to a leaf, along with the
    leaf; on each $2$-weighted edge which does not lead to a leaf, 
    simply remove the weight.
\item The only weights now remaining are even weights $>2$. If such a
  weight is on an edge to a leaf, just remove the weight. If it is on
  an edge joining two nodes, remove the weight and add a white vertex
  in the middle of the edge.
  \end{enumerate}
\end{theorem}
\begin{proof}
  By Theorem \ref{graph links} our graph manifold $G$ is a link
  complement. Eisenbud and Neumann in \cite{EisenbudNeumann} classify
  link complements (in arbitrary homology spheres) in terms of
  ``splice diagrams.'' We first recall from \cite{EisenbudNeumann}
  how to write down the splice diagram in our special case. The splice
  diagram for the $(2,n)$--torus link, in which arrowheads correspond
  to components of the link, is as follows:
  \begin{align*}
    \xymatrix@R=12pt@C=36pt@M=0pt@W=0pt@H=0pt{\ar@{<->}[rr]&&&n=2\\ 
&\Circ\lineto[dd]^(.65)k\\ &&&n=2k>2\\
&\Circ\ar@{->}[l]\ar@{->}[r]&\\
\Circ\lineto[r]^(.65)2&\Circ\lineto[r]^(.35)n\ar@{->}[dd]&\Circ\\
 &&&n=2k+1>2\\&
}
  \end{align*}
  (Omitted splice diagram weights are $1$.) The splice diagram for a
  connected sum of two links is obtained by joining the splice
  diagrams for each link at the arrowheads corresponding to the link
  components along which connected sum is performed, changing the
  merged arrowhead into an ordinary vertex, and adding a new
  $0$--weighted arrow at that vertex.  For example the splice diagram
  corresponding to the Artin presentation graph
$$
  \xymatrix@R=6pt@C=36pt@M=0pt@W=0pt@H=0pt{
\Circ\lineto[r]^2&\Circ\lineto[r]^{4}&\Circ\lineto[r]^3&\Circ}$$ 
would be
$$
\xymatrix@R=24pt@C=36pt@M=0pt@W=0pt@H=0pt{&&\Circ&&\Circ\\
&\Circ\ar@{->}[l]\ar@{->}[u]^(.25)0\lineto[r]&\Circ\lineto[u]^(.25)2\lineto[r]&
\Circ\ar@{->}[u]^(.25)0\lineto[r]&\Circ\lineto[u]_(.25)2\lineto[d]^(.25)3\\
&&&&\Circ
}
$$
Now the nodes of the splice diagram correspond to Seifert pieces in
the geometric decomposition of the graph manifold. Thus the colored
decomposition graph is obtained by taking the full subtree on the nodes of
the diagram with nodes that had arrowheads
attached colored black and the others colored white. This is as
described in the theorem.
\end{proof}


\begin{thebibliography}{BDM}

\bibitem
{BehrstockDrutuMosher:thick}
J.~Behrstock, C.~Dru\c{t}u, and L.~Mosher.
\newblock {Thick metric spaces, relative hyperbolicity, and quasi-isometric
  rigidity}.
\newblock Preprint, \textsc{arXiv:math.GT/0512592}, 2005.

\bibitem
{Brady+:Rigidity}
Noel Brady, Jonathan P.  McCammond, Bernhard M\"uhlherr, Walter D. Neumann.
Rigidity of Coxeter groups and Artin groups. 
In {\em Proceedings of the Conference on Geometric and Combinatorial Group Theory, Part I (Haifa, 2000)}.
Geom. Dedicata {\bf 94} (2002), 91--109.


\bibitem
{Brunner}
A.~M.~Brunner, 
Geometric quotients of link groups.
{Topology Appl.} {\bf 48} (1992), 245--262.

\bibitem{CannonCooper}
J.~Cannon and D.~Cooper.
\newblock A characterization of cocompact hyperbolic and finite-volume
  hyperbolic groups in dimension three.
\newblock Trans. AMS, {\bf 330} (1992), 419--431.


\bibitem{EisenbudNeumann}
D. Eisenbud and Walter D. Neumann.
{\em Three-Dimensional Link Theory and
Invariants of Plane Curve Singularities}. 
Annals of Math. Studies {\bf 110} (Princeton Univ. Press 1985).

\bibitem{EskinFisherWhyte}
Alex~Eskin, David~Fisher, and Kevin~Whyte.
\newblock Quasi-isometries and rigidity of solvable groups.
\newblock Preprint, \textsc{arXiv:math.GR/0511647}, 2005.


\bibitem{Gersten:divergence3mflds}
S. M. Gersten.
\newblock {Divergence in {$3$}-manifold groups}.
\newblock {Geom. Funct. Anal.} {\bf 4} (1994), 633--647.


\bibitem{Gordon:Coherence}
C. McA. Gordon.
\newblock {Artin groups, 3-manifolds and coherence}.
\newblock(Bol. Soc. Mat. Mexicana (3) {\bf10} (2004), Special Issue
in Honor of Francisco "Fico" Gonzalez-Acuna, 193--198.
  
\bibitem{Gromov:PolynomialGrowth}
  M.~Gromov.
  \newblock Groups of polynomial growth and expanding maps.
  \newblock IHES Sci. Publ. Math., {\bf 53} (1981) 53--73.

\bibitem{Gromov:Asymptotic}
M.~Gromov.
\newblock {Asymptotic invariants of infinite groups}.
\newblock In {\em Geometric Group Theory, Vol.
  2 (Sussex, 1991)} {\bf 182} {\em LMS Lecture Notes}, 
  (Cambridge Univ. Press, 1993) 1--295.

\bibitem
{HermillerMeier}
Susan Hermiller and John Meier.
Artin groups, rewriting systems and three-manifolds.
{J. Pure Appl. Algebra.} {\bf 136} (1999), 141--156.

\bibitem
{KapovichLeeb:haken}
M.~Kapovich and B.~Leeb.
\newblock {Quasi-isometries preserve the geometric decomposition of {H}aken
  manifolds}.
\newblock {Invent. Math.} {\bf 128} (1997), 393--416.

\bibitem
{KapovichLeeb:3manifolds}
M.~Kapovich and B.~Leeb.
\newblock {{$3$}-manifold groups and nonpositive curvature}.
\newblock {Geom. Funct. Anal.} {\bf 8} (1998), 841--852.

\bibitem{mckay}Brendan McKay, ``nauty'': a program for isomorphism and
  automorphism of graphs, \verb+http://cs.anu.edu.au/~bdm/nauty/+

\bibitem{milnor} J. Milnor.
A note on curvature and the fundamental group. J. Diff. Geom. {\bf 2}
(1968), 1--7.

\bibitem
{Neumann:commens}
Walter D.~Neumann.
\newblock {Commensurability and virtual fibration for graph manifolds}.
\newblock {Topology} {\bf 36} (1997), 355--378.

\bibitem{neumann-swarup}
Walter D. Neumann and G.A. Swarup. Canonical decompositions of 3-manifolds.
Geometry and Topology {\bf1} (1997), 21--40.

\bibitem{PapasogluWhyte:ends}
P. Papasoglu and K. Whyte.
\newblock {Quasi-isometries between groups with infinitely many ends}.
\newblock {\em Comment. Math. Helv.}, {\bf 77 (1)} (2002) 133--144.


\bibitem
{Perelman:Geom1}
G. Perelman.
\newblock {The entropy formula for the Ricci flow and its geometric 
applications}.
\newblock Preprint, \textsc{arXiv:math.DG/0211159}, 2002.

\bibitem
{Perelman:Geom2}
G. Perelman.
\newblock {Ricci flow with surgery on three-manifolds}.
\newblock Preprint, \textsc{arXiv:math.DG/0303109}, 2003.

\bibitem
{Perelman:Geom3}
G. Perelman.
\newblock {Finite extinction time for the solutions to the Ricci flow on 
certain three-manifolds}.
\newblock Preprint, \textsc{arXiv:math.DG/0307245}, 2003.

\bibitem{Rieffel:H2crossR}
E.~Rieffel.
\newblock {Groups quasi-isometric to $\mathbf {H}^{2}\times\mathbf {R}$}.
\newblock {\em J. London Math. Soc. (2)}, {\bf 64 (1)} (2001) 44--60.

\bibitem{Schwartz:RankOne}
R.~Schwartz.
\newblock The quasi-isometry classification of rank one lattices.
\newblock {\em IHES Sci. Publ. Math.}, {\bf 82} (1996) 133--168.

\bibitem
{Scott:Coherent}
G.~P. Scott.
\newblock {Finitely generated {$3$}-manifold groups are finitely presented}.
\newblock {J. London Math. Soc. (2)} {\bf 6} (1973), 437--440.

\bibitem{svarc} A.S. \v Svarc.
Volume invariants of coverings. Dokl. Akad. Nauk. SSSR {\bf 105}
(1955), 32--34.

\end{thebibliography}
\end{document}